\documentclass[letterpaper,conference]{ieeeconf}  

\IEEEoverridecommandlockouts                              
\usepackage{comment}
\usepackage{graphicx,color}
\usepackage{amsmath}
\usepackage{amssymb}

\newtheorem{eg}{Example}[section]
\newtheorem{definition}[eg]{Definition}
\newtheorem{thm}[eg]{Theorem}

\newtheorem{remark}[eg]{Remark}
\newtheorem{prop}[eg]{Proposition}

\newtheorem{assumption}[eg]{Assumption}
\newenvironment{pf}{\noindent {\bf Proof:}}{\hspace*{\fill}$\square$}

\newtheorem{defn}[eg]{Definition}

\newtheorem{cor}[eg]{Corollary}

\renewcommand{\L}{\mathcal{L}}
\newcommand{\X}{\mathcal{X}}
\newcommand{\D}{{\mathcal D}}
\newcommand{\Z}{\mathcal{Z}}
\newcommand{\Y}{\mathcal{Y}}
\newcommand{\A}{A}
\newcommand{\E}{E}
\renewcommand{\forall}{{\rm \; for \; all \; }}

\newcommand{\la}{\langle}
\newcommand{\ra}{\rangle}

\newcommand{\re}{{\,\rm Re} \,  }
\newcommand{\ran}{{\,\rm ran} \,  }

\begin{document}
\title{On solvability of dissipative partial differential-algebraic equations
}
\author{Birgit Jacob and Kirsten Morris
\thanks{Dept. of Mathematics, Univ. of Wuppertal, Wuppertal, Germany (BJ), Dept. of Applied Mathematics (KM), Univ. of Waterloo, Waterloo, ON, Canada 
         { \tt\small kmorris@uwaterloo.ca}  }%
\thanks{Financial  support of  Natural Sciences and Engineering Research Council of Canada (NSERC) for this research is gratefully acknowledged.}
}
\maketitle

\begin{abstract}
In this article we investigate the solvability of infinite-dimensional differential algebraic equations. Such equations often arise as partial differential-algebraic equations (PDAEs).  A decomposition of the state-space that leads to an extension of the Hille-Yosida Theorem on Hilbert spaces for these equations is described. For dissipative partial differential equations the famous Lumer-Phillips  generation theorem characterizes solvability and also boundedness of the associated semigroup. 
An extension of  the Lumer-Phillips  generation theorem to dissipative differential-algebraic equations is given. The   results is illustrated by coupled systems and  the Dzektser equation.
\end{abstract}

\section{Introduction}

We consider infinite-dimensional differential-algebraic equations (DAEs) 
\begin{align}\label{eqn1}
\frac{d}{dt}Ex(t) & = Ax(t), \quad t\ge 0, \qquad Ex(0)=z_0.
\end{align}
Here $A$ and $E$ are linear operators from $\X$ to $\Z$ and $z_0\in\Z$, where $\X$ and $\Z$ are  complex Hilbert spaces. The operator $E$ is bounded from $\X,$ but $A$ is densely defined and closed on $\X. $
Such equations arise from the coupling of partial differential equations where one sub-system is in equilibrium, as well as in some other situations.  

Establishing well-posedness of these equations, particularly when $E$ is not invertible, is non-trivial.
Solvability of infinite-dimensional differential-algebraic equations has been intensively studied, see \cite{Reis2008, Tro20,ThaTha01,ThaTha96,FavYag04,Sho10, Yag91,FavYag99}. 
For example, Trostorff \cite{Tro20} and Reis and Tischendorf \cite{ReisTischendorf} provide sufficient condition in terms of Hille-Yosida type resolvent estimates.
  In \cite{ThaTha96} the splitting $\X=\ker E\oplus \overline{\ran E^*}$ and $\Z=\ker E^*\oplus \overline{\ran E}$, and the restriction of equation \eqref{eqn1} to this splitting and the solvability is investigated.

In \cite{SFbook} a concept called $(E,r)$-radiality is introduced that also leads to Hille-Yosida type conditions for generation of a semigroup. Under associated conditions, there exists a splitting of $\X$ and $\Z$ into   the kernel of $E$ and a non-orthogonal complement. 
With respect to this splitting the differential-algebraic equations \eqref{eqn1} are written as
\begin{align}\label{Kron}
 \frac{d}{dt}\begin{bmatrix}   0 & 0 \\ 0& E_1 \end{bmatrix}  \begin{bmatrix}   x^0  \\   x^1 \end{bmatrix}   &=  \begin{bmatrix}  A_0 & 0 \\ 0 &A_1 \end{bmatrix} \begin{bmatrix}  x^0 \\ x^1 \end{bmatrix},  \quad t\ge 0,\\ \qquad E_1x^1(0)&=(x_0)^1. \nonumber
\end{align}
where $E_1:\X^1\rightarrow \Z^1$ is bounded and invertible, $A_0:\D(A)\cap \X^0\rightarrow \Z^0$ is closed and invertible, $A_1:\D(A)\cap \X^1\rightarrow \Z^1$ is closed, and $A_1 E_1^{-1}$ generates a $C_0$-semigroup in $\Z^1$. Here $\D(A)$ denotes the  domain of the operator $A$.

 As we work with Hilbert spaces instead on general Banach spaces, we are able to simplify and weaken the required conditions considerably. 

The difficulty with this approach, as with the classical Hille-Yosida Theorem, and with other resolvant estimates, is that it can be difficult to confirm that the assumptions are satisfied. 
The Lumer-Phillips Theorem e.g. \cite{CZbook} is a very useful tool  in the standard $E=I$ situation for establishing that an operator $A$ generates a $C_0$-semigroup. 
Favini and Yagi \cite[Page 37]{FavYag99} show that if
    $(\lambda E-A)^{-1}\in \L(\Z,\X)$ for some $\lambda >0$ and 
$\re\langle Ax, Ex\rangle_\Z\le 0$ for all $ x\in \D(A)$,
 then  for every $x_0\in E(\D(A))$ the DAE \eqref{eqn1}  has a unique classical solution; that is  $x:[0,T]\rightarrow \D(A)$  and  $Ex\in C^1([0,T];\Z)$,  $Ax\in C([0,T];Z)$,  and DAE \eqref{eqn1} is satisfied.  They further provide results for parabolic DAEs. However, they do not investigate a splitting of the state space as in \eqref{Kron}, nor do they show generation of a $C_0$-semigroup. 
 The main result of this paper is a generalization of the Lumer-Phillips Theorem  to dissipative infinite-dimensional DAEs.

The framework  of $(E,r)$-radiality is summarized and adapted to the Hilbert space situation in the next section. Some new results are proven.   In Section 3 we prove a Lumer-Phillips Theorem  for dissipative infinite-dimensional DAEs. Well-posedness of a class of  coupled systems  is shown in Section 4.  Finally in Section 5, the results are applied to the  Dzektser equation.

The following notation will be used. By $\L(\X,\Y)$ we denote the set of all linear and bounded operators form the Hilbert space $\X$ to the Hilbert space $\Y$. We denote the kernel of an operator $A$ by $\ker A$ and its range by $\ran A$. If an operator $A$ is closable, the we denote the closure by $\overline{A}$.

\section{Radiality and semigroup generation}


Let $\X$, $\Z$ be  complex Hilbert spaces, $E\in \L(\X,\Z)$, $A:D(A)\subset \X\rightarrow \Z$  densely defined and closed, and $z_0\in \Z$. 
By
$$
\varrho(E,A) :=\{ s\in \mathbb C\mid (sE-A)^{-1}\in \L(\Z,\X)\}
$$
we denote the   {\it{resolvent set}} of the operator pencil $(E,A)$. 
For $s\in  \varrho(E,A)$, we define the right- and left-$E$ resolvents of $A$  (with respect to $E$) by
$$ R^E (s,A) = (sE-A)^{-1} E , \quad  L^E (s,A) = E (sE-A)^{-1},$$
In particular,  if $0 \in \varrho(E,A), $ $$R^E (0,A) = -A^{-1} E, \quad L^E(0,A) = -E A^{-1}.$$
\begin{defn}
The operator $A$ is {\em $E$-radial} if 
\begin{itemize}
\item
$s \in \varrho (E,A)$ for all real $s>0$,
\item
there exists $K>0$ such that for all $n\in \mathbb N$ and for all real $s>0$
\begin{align} 
 \| \left( R^E (s,A) \right)^n \|_{\L(\X,\X)} &\leq \frac{K}{s^n}, \label{eq-rada} \\
 \| \left( L^E (s,A) \right)^n  \|_{\L(\Z,\Z)}&\leq \frac{K}{s^n},\label{eq-radb}
\end{align}
\end{itemize}
\end{defn}
\vspace{1ex}

In the case $E=I$ statements \eqref{eq-rada} and \eqref{eq-radb} are equivalent; and either statement implies generation of a $C_0$-semigroup by the Hille-Yosida Theorem. 
 In \cite{SFbook} a more general concept $(E,p)$ radiality is considered. In that framework,  an $E$-radial operator is  $(E,0)$-radial. Here, only $p=0$ is considered and also the spaces are assumed to be Hilbert spaces. 

\begin{defn}
The operator $A$ is {\em weakly $E$-radial} if 
$s \in \varrho (E,A)$ for all $s>0$,
and \eqref{eq-rada}-\eqref{eq-radb} holds with $n=1$.
\end{defn}

Clearly any $E$-radial operator is weakly $E$-radial. The converse holds if $K\leq 1.$

\begin{defn}
The operator $A$ is {\em strongly $E$-radial} if it is $E$-radial and there is a linear and dense subspace $\tilde \Z $  of $\Z$ such that 
\begin{equation*}
\begin{aligned} 
   \|  R^E (s,A) (\lambda E- A )^{-1} A x \|_{\X} &\leq \frac{{\rm const} (x) }{\lambda s } , \quad  x \in D(A) , \\
     \|  A (\lambda E - A )^{-1}   L^E (s,A)z  \|_{\Z}  & \leq   \frac{{\rm const} (z )}{\lambda s}   , \quad \quad   z \in \tilde \Z. 
 \end{aligned}
  \label{eq-srad} \end{equation*}
\end{defn}
\vspace{1ex}

Define for some $\alpha \in \varrho (E,A)$, 
\begin{align*}
\X^0 &= \ker R^E (\alpha,A)=\ker E,\\
 \Z^0 &=  \ker L^E (\alpha,A)=\{Ax\mid x\in D(A)\cap \ker E\},\\
\X^1 & = \overline{\ran R^E (\alpha,A) },\\
 \Z^1 &= \overline{ \ran  L^E (\alpha,A)}.
\end{align*}
These spaces are independent of the choice of $\alpha$ (\cite[Lem.~2.1.2, pg.~18]{SFbook}) and  also  if $A$ is weakly $E$-radial, then
\begin{align*}
\lim_{s \to \infty} s R^E (s,A) x  = x , \quad \forall x \in \X^1,\\
\lim_{s \to \infty} s L^E (s,A) z  = z , \quad \forall z \in \Z^1,
\end{align*}
see \cite[Lem. 2.2.6]{SFbook}.
If $A$ is weakly $E$-radial, then since the Hilbert spaces $\X$ and $\Z$ are reflexive,  \cite[Theorem 2.5.1]{SFbook}
implies
$$ \X=\X^0\oplus \X^1 \qquad \text{  and } \qquad \Z=\Z^0\oplus \Z^1.$$ 
Thus if $A$ is weakly $E$-radial,
\begin{itemize}
\item
$ P : \X \to \X$ defined by $$Px:=\lim_{s \to \infty} s R^E (s,A) x$$ is  a projection onto $\X^1$ with $\ker P = \X^0$ and $\ran P=\X^1$, 
\item
$ Q: \X \to \X$ defined by $$Qz:=\lim_{s \to \infty} s L^E (s,A) z$$ is  a projection onto $\Z^1$ with $\ker Q = \Z^0$ and $\ran Q=\Z^1$.
\end{itemize}
In general, both $P$ and $Q$ are non-orthogonal projections.

Define restrictions of $E$ and $A$ as follows:
$$E_0 := E\vert_{\X^0} , \quad A_0 := A\vert_{D (A_0)  }, \; D (A_0) = \X^0 \cap D(A), $$
$$E_1 := E\vert_{\X^1} , \quad A_1 := A\vert_{D (A_1)  }, \; D (A_1) = \X^1 \cap D(A). $$
In \cite[Lem.~2.2.1, pg.~20]{SFbook} it is shown that 
 $E_0 \in \L  (\X^0 , \Z^0 )$ and  $ A_0: D (A_0)  \to \Z^0 $. 
Further, if $A$ is weakly $E$-radial, then $A_0$ is boundedly invertible;  that is,
$$A_0^{-1} \in \L (\Z^0,\X^0),$$ 
see \cite[Lem. 2.2.4, pg.~22]{SFbook} and also
$$A_0^{-1} E_0=0, \quad E_0 A_0^{-1}=0$$
on $\X^0$ and $\Z^0$, respectively, by  \cite[Lem.~2.2.5, pg.~22]{SFbook}.

The following proposition has been proved in \cite[Cor. 2.5.1,pg 38]{SFbook} if $A$ is strongly radial. We are able to weaken this assumption because we deal with Hilbert spaces.
\begin{prop}\label{thm-QP}
If A is weakly $E$-radial, then
\begin{enumerate}
\item for all $x \in D (A) , $  $P x \in D (A) $ and $AP x = QA x ,$
\item
for all $x \in \X, $ $EP x = Q E x .$
\end{enumerate}
\end{prop}

\begin{pf}
Recall that the operator $P$ is defined by
$$Px=\lim_{s \to \infty} s R^E (s,A) x.$$
For any $x \in D (A) \subset \X$,  by \cite[Equation (2.1.8), pg.~17]{SFbook} 
$$ A R^E (s,A) x = L^E (s,A) A x . $$
Let $x\in D(A)$. Since $R^E(s,A) x \in D (A) , $ and $A$ is closed, $Px \in D (A)$. 
Thus 
\begin{align*}
 A P x &= A (\lim_{s \to \infty} s R^E (s,A) x )\\
&= \lim_{s \to \infty} A s R^E (s,A) x\\
&=\lim_{s \to \infty} s L^E (s,A) A x  \\
&= Q A x. 
 \end{align*}
 This proves Part 1). Part 2) follows easily using the fact that $E \in \L (\X, \Z )$. For any $x \in \X,$
 \begin{align*}
 EP x &= E\lim_{s \to \infty} s R^E (s,A) x\\
 &= \lim_{s \to \infty} sE  R^E (s,A) x\\
&= \lim_{s \to \infty}s  L^E (s,A) E x\\ 
 &= Q E x.
 \end{align*}
This concludes the proof. \end{pf}

Thus, if A is weakly $E$-radial, then the operators $A_0$, $A_1$, $E_0$ and $E_1$ are invariant with repect to the projected spaces. More precisely,  if A is weakly $E$-radial, by \cite[Lem.~2.2.1, pg.~20 and Cor. 2.5.2, pg. 39]{SFbook}
\begin{itemize}
\item $E_0 \in \L (\X^0 , \Z^0)$,
\item $E_1 \in \L (\X^1 , \Z^1)$,
\item $ A_0: D (A_0)\subset \X^0  \to \Z^0$ is densely defined, closed, and boundedly invertible,
\item $ A_1: D (A_1)\subset \X^1  \to \Z^1$ is densely defined and closed.
\end{itemize}

%
The following proposition was proven in \cite[Cor. 2.5.3,pg.~40]{SFbook} with an assumption that $A$ is strongly radial.

\begin{prop}
If $A$ is weakly $E$-radial and $\ran E$ is closed in $\Z$, then $E_1 \in \L (\X^1 , \Z^1)$ is boundedly invertible.
\end{prop}
\begin{pf}
The fact that $E_1 \in \L (\X^1 , \Z^1)$ follows from \cite[Cor. 2.5.2, pg. 39]{SFbook}. Since $\X=\X^0\oplus \X^1$, and $\X^0 = \ker E , $ it follows that $E_1$ is injective. Thus it remains to show that $E_1$ is surjective. By the definition of $\Z_1 ,$
$$ \Z^1=\overline{E(D(A))} \subset \ran E=\ran E_1 . $$
Since $ \ran E_1  \subset \Z^1 ,$
$E_1$ is surjective. Thus $E_1$ has an inverse defined on all of $\Z^1$ and  so by the Closed Graph Theorem this inverse is bounded. \end{pf}

This framework allows us to  decompose the system using the non-orthogonal projections $P$ and $Q$ if  $A$ is weakly $E$-radial and $\ran E$ is closed.
Defining
\begin{align*} 
\tilde P &= \begin{bmatrix} I -P \\ P  \end{bmatrix}  \in \L (\X, \X^0 \times \X^1 ),\\
\tilde Q &= \begin{bmatrix} I- Q \\ Q \end{bmatrix} \in \L ( \Z, \Z^0 \times \Z^1),
\end{align*}
we obtain
\begin{align*}
\tilde P^{-1}  &= \begin{bmatrix} I & I  \end{bmatrix}  \in \L (\X^0 \times \X^1, \X ),\\
\tilde Q^{-1} &= \begin{bmatrix} I & I  \end{bmatrix}  \in \L ( \Z^0 \times \Z^1, \Z),
\end{align*}
where $I$ above indicates the natural injection on the various spaces; the different spaces are not explicitly indicated. 
Let $z\in \Z$ and $\begin{bmatrix}  z_0 \\ z_1 \end{bmatrix}:=\tilde P z$. Then we obtain the chain of equivalences
\begin{align*}
&&\quad E  \frac{d}{dt} z &= A z \\
\Longleftrightarrow \quad &&\tilde Q  E \tilde P^{-1} \begin{bmatrix}  \dot z_0 \\  \dot z_1 \end{bmatrix}  &= \tilde Q  A \tilde P^{-1} \begin{bmatrix}  z_0 \\ z_1 \end{bmatrix} \\ 
\Longleftrightarrow \quad&& \begin{bmatrix}   E_0 & 0 \\ 0& E_1 \end{bmatrix}  \begin{bmatrix}  \dot z_0 \\  \dot z_1 \end{bmatrix}   &= \begin{bmatrix}  A_0 & 0 \\ 0 &A_1 \end{bmatrix} \begin{bmatrix}  z_0 \\ z_1 \end{bmatrix}\\
\Longleftrightarrow \quad&& \begin{bmatrix}  0 & 0 \\ 0 &I \end{bmatrix}  \begin{bmatrix}  \dot z_0 \\  \dot z_1 \end{bmatrix}   &=  \begin{bmatrix}  I & 0 \\ 0 & E_1^{-1}A_1 \end{bmatrix} \begin{bmatrix}  z_0 \\ z_1 \end{bmatrix}.
\end{align*}

Our main result in this section is as follows.
\begin{thm}\label{thm:semigroup}
If $A-\alpha E $ is $E$-radial and $\ran E$ is closed, then the operator $E_1^{-1} A_1$  with domain $D (A) \cap \X^1$ generates a  $C_0$-semigroup $(S(t))_{t\ge 0}$  on $\X^1$  with bound $Ke ^{\alpha t} .$
\end{thm}

\begin{pf}
First define $\tilde A = (A- \alpha E) . $
By our assumption the operator $E_1^{-1} \tilde A_1$  with domain $D (A) \cap \X^1$ is well-defined, closed and densely defined. The definition of $E$-radiality further implies $(0,\infty)\in \varrho(E_1^{-1} \tilde A_1)$ and there exists $K>0$ such that for all $n\in \mathbb N$ and for all real $s>0$
$$
 \|(sI-E_1^{-1} \tilde A_1)^n \|_{\L(\X^1,\X^1)} \leq \frac{K}{s^n}.
$$
Thus from the Hille-Yosida Theorem, $E_1^{-1} \tilde A_1$ generates a bounded $C_0$-semigroup on $\Z^1$ with bound $K. $ 

Now note   $\tilde A_1 = (A- \alpha E)_1 = A_1 - \alpha E_1 .$   The statement of the theorem now follows. 
\end{pf}

Alternatively,  the same projections can be applied to  
$$  \frac{d}{dt}E z = A z . $$
This yields
\begin{align*}
&&\quad \frac{d}{dt}E z &= A z \\
\Longleftrightarrow \quad && \frac{d}{dt}\tilde Q  E \tilde P^{-1} \begin{bmatrix}   z_0 \\   z_1 \end{bmatrix}  &= \tilde Q  A \tilde P^{-1} \begin{bmatrix}  z_0 \\ z_1 \end{bmatrix} \\ 
\Longleftrightarrow \quad&& \frac{d}{dt}\begin{bmatrix}   E_0 & 0 \\ 0& E_1 \end{bmatrix}  \begin{bmatrix}  z_0 \\   z_1 \end{bmatrix}   &= \begin{bmatrix}  A_0 & 0 \\ 0 &A_1 \end{bmatrix} \begin{bmatrix}  z_0 \\ z_1 \end{bmatrix}\\
\Longleftrightarrow \quad&& \frac{d}{dt} \begin{bmatrix}  0 & 0 \\ 0 &I \end{bmatrix}  \begin{bmatrix}  x_0 \\   x_1 \end{bmatrix}   &=  \begin{bmatrix}  I & 0 \\ 0 & A_1E_1^{-1} \end{bmatrix} \begin{bmatrix}  x_0 \\ x_1 \end{bmatrix},
\end{align*}
where  $\left[\begin{matrix}  z_0 \\ z_1 \end{matrix}\right]=\left[\begin{matrix}  A_0^{-1}x_0 \\ E_1^{-1} x_1 \end{matrix}\right]$.
\vspace{1ex}

\begin{cor}
If $A-\alpha I $ is $E$-radial and $\ran E$ is closed, then the operator $ A_1 E_1^{-1}$  with domain $E_1 (\D(A) \cap \X_1) $  generates a $C_0$ semigroup on $\Z_1 $ with bound $K e^{\alpha t}.$
\end{cor}

Thus well-posedness of the abstract Cauchy problem $ E \frac{d}{dt} x = A x$  reduces to well-posedness of  $  \dot x_1 = E_1^{-1}A_1 x_1$. Further,   well-posedness of  $ \frac{d}{dt}(E x) = A x$ is reduced to well-posedness of  $ \dot z_1 = A_1E_1^{-1} z_1. $

\section{Dissipative pencils}

The fact whether an operator $A$ is  $E$-radial is in general not easy to verify. Therefore, in this section we aim to generalize the Lumer-Phillips theorem to dissipative differential-algebraic equations. We start with the following definition. As in the previous section $\X$ and $\Z$ are  complex Hilbert spaces, $E\in \L(\X,\Z)$, and $A:\D (A)\subset \X\rightarrow \Z$  is densely defined and closed. 

\begin{definition}
The operator pencil $(E,A)$ is called {\rm dissipative}, if
\begin{equation}\label{eqn:diss}
\|(\lambda E-A)x\|_\Z \ge \lambda \|Ex\|_\Z, \qquad \lambda >0, x\in \D(A).
\end{equation}
\end{definition}
In the following three propositions  we define dissipative operator pencils and prove some implications.
\begin{prop}\label{prop:diss}
The following statements are equivalent:
\begin{enumerate}
\item The operator pencil $(E,A)$ is dissipative.
\item  Re$\langle Ax, Ex\rangle_\Z\le 0$, $x\in \D(A)$,
\end{enumerate}
\end{prop}
\begin{pf}
For $\lambda>0$ and $x\in \D(A)$ equation \eqref{eqn:diss} is equivalent to
\begin{equation*}
\frac{1}{2\lambda}   \| A x \|^2  \ge   \re \la Ax , Ex \ra.
\end{equation*}
This implies the statement of the proposition.
\end{pf}

\begin{remark}
If $\frac{d}{dt}Ex(t)  = Ax(t)$, $t\ge 0$,  has a classical solution and satisfies
$
\frac{d}{dt}\|Ex(t)\|^2 \le 0$
for every classical solution $x$, then the operator pencil is dissipative. 
\end{remark}
\begin{prop}\label{prop:ker}
If the  operator pencil $(E,A)$ is dissipative, then 
the following are equivalent
\begin{enumerate}
\item $\ker A \cap \ker E=\{0\}$,
\item $\ker (\lambda E-A)=\{0\}$ for one $\lambda>0$,
\item $\ker (\lambda E-A)=\{0\}$ for every $\lambda>0$,
\end{enumerate}
\end{prop}

\begin{pf}
Clearly statement (3) ~implies (2) ~and (2)~implies (1). Thus it remains to show that (1) implies (3). Assume that (3) does not hold, that is, there exists  $\lambda>0$ and $x\in \ker (\lambda E-A)$, $x\not=0$. This implies 
\begin{equation*}
\lambda^2 \|Ex\|^2 -2\lambda \re \la Ex,Ax\ra + \|Ax\|^2=0.
\end{equation*}
Since $(E,A)$ is dissipative,  all terms of the right hand side are non-negative. Thus $Ex=0$ and $Ax=0$, which implies that (1) does not hold. Hence statement (1) implies (3) . 
\end{pf}

\begin{prop}\label{prop:spectrum}
If the  operator pencil $(E,A)$ is dissipative and $\varrho(E,A)\cap (0,\infty)\not=\emptyset$, then $(0,\infty) \subset\varrho(E,A)$.  
\end{prop}
\begin{pf}
Let $\lambda_0 \in \varrho(E,A)\cap (0,\infty)$. 
Proposition \ref{prop:ker} implies 
$\ker(\lambda E-A)=\{0\}$ for every $\lambda>0$. 

Next we show that $\ran(\lambda E-A)$ is dense in $\Z$ for every $\lambda>0$. Let $\lambda>0$ be arbitrary and $z\in \Z$ be orthogonal to $\ran(\lambda E-A)$. Since $\ran(\lambda_0 E-A)=\Z$ there exists $x\in\D(A)$ such that $z=(\lambda_0 E-A)x$. This implies
\begin{align*}
0 =&\re \la z, (\lambda E-A)x\ra\\
&=\re\la(\lambda_0 E-A)x , (\lambda E-A)x\ra \\
&= \lambda_0\lambda \|Ex\|^2 -(\lambda+\lambda_0) \re \la Ex,Ax\ra + \|Ax\|^2.
\end{align*}
Because  $(E,A)$ is dissipative, all terms of the right hand side are non-negative. Thus $Ex=0$ and $Ax=0$, which implies $\lambda_0Ex-Ax=0$.  Since   $\ker(\lambda_0 E-A)=\{0\}$, it follows that   $x=0$ and therefore $z=0$. Thus $\ran(\lambda E-A)$ is dense in $\Z$ for every $\lambda>0$.

It remains to show that  $\ran(\lambda E-A)$ is closed in $\Z$ for every $\lambda>0$.
We note that for an injective, closed and densely defined operator $T:\D(T)\subset \X\rightarrow \Z$ the following statements are equivalent
\begin{itemize}
\item $\ran (T)$ is closed in $\Z$.
\item There exists $c>0$ such that $\|Tx\|\ge c\|x\|$ for every $x\in \D(T)$.
\end{itemize}
Thus there exists $c_{\lambda_0}>0$ such that $\|(\lambda_0 E-A)x\|\ge c_{\lambda_0}\|x\|$ for every $x\in \D(A)$. This implies that  for any $x\in \D(A)$ with $\|x\|=1$:
\begin{equation}\label{eqn:closed}
\lambda_0^2 \|Ex\|^2 -2\lambda_0 \re \la Ex,Ax\ra + \|Ax\|^2\ge c_{\lambda_0}.
\end{equation}

It remains to show that for every $\lambda>0$  there exists $c_{\lambda}>0$ such that
\begin{equation}
\lambda^2 \|Ex\|^2 -2\lambda \re \la Ex,Ax\ra + \|Ax\|^2\ge c_\lambda \label{eq-cl}
\end{equation}
for all $x\in \D(A)$ with $\|x\|=1$. Assume that this is not true: that is, there exists $\lambda>0$ and a sequence $(x_n)_n\subset \D(A)$ with $\|x_n\|=1$ such that 
\begin{equation*}
\lambda^2 \|Ex_n\|^2 -2\lambda \re \la Ex_n,Ax_n\ra + \|Ax_n\|^2 \rightarrow 0
\end{equation*}
as $n\rightarrow \infty$. As $(E,A)$ is dissipative, all three terms are non-negative and thus $ \|Ex_n\|^2\rightarrow 0$,  $ \|Ax_n\|^2\rightarrow 0$ and 
 $ \la Ex_n,Ax_n\ra\rightarrow 0$ as $n\rightarrow \infty$, which implies that  \eqref{eqn:closed} does not hold. Thus, statement \eqref{eq-cl} holds. 
\end{pf}

%
The first main result of this section is summarized in the following theorem.

\begin{thm}
\label{thm-Ediss}
 If    $\lambda \in \varrho(E,A)$ for some $\lambda >0$ and 
 \begin{align*}
 \re\langle Ax, Ex\rangle_\Z&\le 0, \quad x\in \D(A),\\
  \re\langle A^*x, E^*x\rangle_\X&\le 0, \quad x\in \D(A^*),
 \end{align*}
 then 
\begin{enumerate}
\item $(0,\infty)\subset \varrho(E,A)$ and $(0,\infty)\subset \varrho(E^*,A^*)$.
\item $\|E(\lambda E-A)^{-1}\| \le \frac{1}{\lambda}$ for  $\lambda>0$,
\item $\|(\lambda E-A)^{-1}E\| \le \frac{1}{\lambda}$ for  $\lambda>0$,
\item $A$ is $E$-radial.
\end{enumerate}
\end{thm}

\begin{pf}
The first statement follows from Proposition \ref{prop:spectrum}, and the second  from Proposition \ref{prop:diss}. Moreover, Proposition \ref{prop:diss} implies 
\begin{equation*}
\|E^*(\lambda E^*-A^*)^{-1}\| \le \frac{1}{\lambda}
\end{equation*}
 for  $\lambda>0$, which implies the third statement. The last statement now follows directly from the definition of $E$-radiality.
%
\end{pf}

Theorem \ref{thm-Ediss} together with  Theorem \ref{thm:semigroup}  now implies the second result of this section. 

\begin{thm}\label{thm:Erad}
If  the operator $E$ has closed range,   $\lambda \in \varrho(E,A)$ for some $\lambda >0$ and 
 \begin{align*}
 \re\langle Ax, Ex\rangle_\Z&\le 0, \quad x\in \D(A),\\
  \re\langle A^*x, E^*x\rangle_\X&\le 0, \quad x\in \D(A^*),
 \end{align*}
Then the Hilbert spaces $\X$, $\Z$  can be split as
$ \X=\X^0\oplus \X^1$   and $\Z=\Z^0\oplus \Z^1$,
where
\begin{align*}
\X^0 &= \ker E,&
 \Z^0 &= \{Ax\mid x\in \D(A)\cap \ker E\},\\
\X^1 & = \overline{\ran R^E (\alpha,A) },&
 \Z^1 &=  \ran  E
\end{align*}
for some (and hence every) $\alpha \in \rho (E, A)$. 
Further,
$ P : \X \to \X$ defined by $Px:=\lim_{s \to \infty} s R^E (s,A) x$ is  a projection onto $\X^1$ with $\ker P = \X^0$, and
$ Q: \X \to \X$ defined by $Qz:=\lim_{s \to \infty} s L^E (s,A) z$ is  a projection onto $\Z^1$ with $\ker Q = \Z^0$. For 
\begin{align*}
E_1 &:= E\vert_{\X^1}: \X^1\rightarrow \Z^1,\\
A_1 &:= A\vert_{\D (A_1)}:\D (A_1)\subset \X^1\rightarrow \Z^1, \\
 \D (A_1) &= \X^1 \cap \D(A),
\end{align*}
we have $E_1$ is boundedly invertible, 
and $A_1$ is closed and densely defined. 

The equation $\frac{d}{dt}(E  x) = A x$ is equivalent to
\begin{align*}
 \frac{d}{dt}\begin{bmatrix}   0 & 0 \\ 0& I \end{bmatrix}  \begin{bmatrix}   x^0  \\   x^1 \end{bmatrix}   &=  \begin{bmatrix}  I & 0 \\ 0 & E_1^{-1} A_1 \end{bmatrix} \begin{bmatrix}  x^0 \\ x^1 \end{bmatrix},
\end{align*}
where $x= \left[\begin{smallmatrix}  x^0 \\ x^1 \end{smallmatrix}\right]$, and $ E_1^{-1} A_1 $ generates a contraction  semigroup on $\X^1$.
\end{thm}

\section{Coupled systems}

In this section we study an application of  Theorem \ref{thm:semigroup}  to a class of  differential-algebraic systems defined on  a  Hilbert space $\Z ,$
\begin{align*}
\frac{d}{dt} x(t) &= A_1 x (t) + A_2 y (t) \\
0 &= A_3 x (t) + A_4 y (t) ,
\end{align*}
where for $i=1\ldots 4, $ $A_i:\D(A_i)\subset Z\rightarrow Z$  are closed and densely defined.
This class of systems is of the form \eqref{eqn1}:  
\begin{equation} \label{eqn:coupled}
\frac{d}{dt}(\underbrace{ \left[\begin{smallmatrix} I & 0 \\ 0& 0 \end{smallmatrix}\right]}_ {\Large E}x(t)) = \underbrace{\left[\begin{smallmatrix} A_1 & A_2 \\ A_3 & A_4 \end{smallmatrix}\right]}_{\Large A}x(t), \qquad t>0,
\end{equation}
with  $\mathcal X=\mathcal{Z}=Z\times Z$ and 
$$  \D({A})=({\D}(A_1)\cap\D(A_3))\times(\D(A_2)\cap\D(A_4)).$$
Although this is a very particular class of systems, a number of applications fit this class.  For examples see \cite{MO2013} for   piezo-electric beams with  quasi-static magnetic effects, \cite{Chaturvedi}  for lithium-ion cell models and   \cite{Hanke2007} for a model with convection-diffusion dynamics.

In order for $\varrho (E,A)$ to be non-empty, it is necessary that $A_4$ have a bounded inverse. This assumption will be made throughout this section, as well as several other assumptions that will guarantee well-posedness.
\begin{assumption}
\label{asn41}
\begin{description}
\item{(a)}
Let $0\in \varrho(A_4)$,  $\D(A_4)\subset \D(A_2)$ and $\D(A_4^*)\subset \D(A_3^*)$. By \cite[Rem.~2.2.315]{tretter}, this implies that the operator $A_2 A_4^{-1}A_3:\D(A_3)\rightarrow Z$ is well defined. Assume also $\overline{A_2A_4^{-1}A_3}\in \L(Z)$. 
\item{(b)}
Let there exist $M\ge 1$ and $\omega\in \mathbb R$ such that for every $s>\omega, $ $ s\in \varrho(A_1)$ and
\begin{align*}
\|(s-A_1)^{-n}\|\le \frac{M}{(s-\omega)^n}, \quad s>\omega, n\in \mathbb N.
\end{align*}
\end{description}
\end{assumption}
\vspace{1ex}

\begin{remark}
If  $A_2, A_3 \in \L (Z)$, then Assumption \ref{asn41}a reduces to $0\in \varrho(A_4).$ 
\end{remark}
\begin{remark}
By the Hille-Yosida Theorem, Assumption \ref{asn41}b  is equivalent to assuming that $A_1$ is closed and generates a $C_0$-semigroup on $Z$.
\end{remark}

For $\mu>\omega$ 
define the 
\textit{Schur complement } $$S_1(\mu):\D(A_1)\subset Z\rightarrow Z$$ by
\begin{align*}
S_1(\mu)&:=\mu-A_1+\overline{A_2A_4^{-1}A_3}.
\end{align*}
Because  $A_1$ is closed and densely defined, the Schur complement $S_1(\mu)$ is closed and densely defined. Moreover, 
we can factor $S_1(\mu)$ as
\begin{align}
S_1(\mu)&:=(\mu-A_1)[I+(\mu-A_1)^{-1}\overline{A_2A_4^{-1}A_3}]. \label{eq-S1}
\end{align}
A Neumann series argument   yields that  the operator $S_1(\mu)$ is invertible for $\mu>\omega_0:=\omega+ M \|\overline{A_2A_4^{-1}A_3}\|$ and also
\begin{align}
\|S_1(\mu)^{-n}\| &\le \|(\mu-A_1)^{-n}\| \|[I+(\mu-A_1)^{-1}\overline{A_2A_4^{-1}A_3}]^{-1}\|^n\nonumber \\
& \le  \frac{M}{(\mu-\omega)^n}\frac{1}{\left(1- \frac{M \|\overline{A_2A_4^{-1}A_3}\|}{\mu -\omega}\right)^n}\label{schur}\\
&= \frac{M}{\left(\mu-\omega_0\right)^n}.\nonumber
\end{align}
Application of \cite[Thm.~2.3.3]{tretter} yields immediately the statement of the following proposition.

\begin{prop}\label{thm:tretter1}
The operator $\A$ is closable and for every $\mu>\omega_0$ we have that $\mu\in \varrho(E, \overline{A)}$ and
\begin{align*}
(\mu E&- \overline{A})^{-1} =\begin{pmatrix} \mu-A _1 & -A_2 \\ -A_3 & -A_4  \end{pmatrix}^{-1}\\
&=\begin{pmatrix} I & 0 \\ -\overline{A_4^{-1}A_3} & I  \end{pmatrix}\begin{pmatrix} S_1(\mu)^{-1} & 0 \\ 0 & -A_4^{-1}  \end{pmatrix}\begin{pmatrix} I & -A_2A_4^{-1} \\ 0 & I \end{pmatrix}\\
&=\begin{pmatrix} S_1(\mu)^{-1} & -S_1(\mu)^{-1}A_2A_4^{-1} \\ -\overline{A_4^{-1}A_3}S_1(\mu)^{-1} & *  \end{pmatrix}.
\end{align*}
\end{prop}


Using Proposition \ref{thm:tretter1}, 
\begin{align}
((\mu E- A)^{-1}E)^n&=\begin{pmatrix} S_1(\mu)^{-n} & 0 \\ -\overline{A_4^{-1}A_3}S_1(\mu)^{-n} & 0  \end{pmatrix} , \label{eq-R}\\
&=\begin{pmatrix} I & 0 \\ -\overline{A_4^{-1}A_3} & 0  \end{pmatrix} \begin{pmatrix} S_1(\mu)^{-n} & 0 \\ 0 & 0  \end{pmatrix} \nonumber \\
(E(\mu E- A)^{-1})^n&=\begin{pmatrix} S_1(\mu)^{-n} & -S_1(\mu)^{-n}A_2A_4^{-1} \\ 0 & 0  \end{pmatrix} \label{eq-L}\\
&=\begin{pmatrix} S_1(\mu)^{-n} & 0 \\ 0 & 0  \end{pmatrix}\begin{pmatrix} I & -A_2A_4^{-1} \\ 0 & 0  \end{pmatrix}. \nonumber
\end{align}
These calculations together with \eqref{schur} show that $A-\omega_0 E$ is $E$-radial.  Further,  $\ran \E$ is closed. Thus  Theorem \ref{thm:semigroup} is applicable. 

The projections $P$ and $Q$ will be explicitly calculated for this class of systems.
Note that because of Assumption \ref{asn41}, $A_1$ generates a $C_0$-semigroup, and thus 
$$\lim_{s \to \infty}  s (S_1 (s ))^{-1} z = z. $$
This implies
\begin{align*}
P\begin{pmatrix} x\\y \end{pmatrix} &=\lim_{s \to \infty} s (sE- \overline{A})^{-1} \begin{pmatrix} x\\0 \end{pmatrix}\\
&= \lim_{s \to \infty}\begin{pmatrix} sS_1 (s )^{-1} x \\ -\overline{A_4^{-1} A_3}sS_1 (s )^{-1} x \end{pmatrix}\\
&= \begin{pmatrix} I & 0 \\ -\overline{A_4^{-1}A_3} & 0 \end{pmatrix}\begin{pmatrix} x\\y \end{pmatrix}
\end{align*}
and similarly
\begin{align*}
Q\begin{pmatrix} x\\y \end{pmatrix} &= \begin{pmatrix} I & -A_2 A_4^{-1} \\ 0  & 0 \end{pmatrix}\begin{pmatrix} x\\y \end{pmatrix}.
\end{align*}

\section{Example: Dzektser equation}

We consider the Dzektser equation
\begin{align*}
\frac{\partial}{\partial t}\left(1+\frac{\partial^2}{\partial \zeta^2}\right)x(\zeta,t) &=\left(\frac{\partial^2}{\partial \zeta^2}+2\frac{\partial^4}{\partial \zeta^4}\right)x(\zeta,t),
\end{align*}
$t>0$ and $\zeta\in(0,\pi)$, with boundary conditions
\begin{align*}
x(0,t)&=x(\pi,t)=0,\quad t>0\\
 \frac{\partial^2x}{\partial \zeta^2}(0,t)&=\frac{\partial^2 x}{\partial \zeta^2}(\pi,t)=0, \quad t>0.
\end{align*}
Let $\Z=L^2(0,\pi)$ and $\X=H^2(0,\pi)\cap H^1_0(0,\pi)$ with $\|x\|^2_{\X}=\|x''\|^2_{\Z}$, $E\in\L(\X,\Z)$ and $A:\D(A)\subset \X\rightarrow \Z$ given by
\begin{align*}
Ex&=x+x'',\\
Ax&=x''+2x^{(4)},\\
\D(A) &=\{x\in H^4(0,\pi)\cap H^1_0(0,\pi)\mid x''(0)=x''(\pi)=0\}.
\end{align*}
This system was shown in \cite{GeGeZhang} to have a splitting  \eqref{Kron} using the eigenfunction expansion. Here the system will be shown to be dissipative $E-$radial.
For $x\in \D(A)$ we calculate 
\begin{align*}
\re\langle Ax, Ex\rangle_\Z&= \re\int_0^\pi (x''+2x^{(4)})(\overline{x}+\overline{x}'')d\zeta\\
&= -\|x'\|^2_{L^2(0,\pi)}+\|x''\|^2_{L^2(0,\pi)}\\
&\quad  -2\|x^{(3)}\|^2_{L^2(0,\pi)}-2\re\int_0^\pi x^{(3)}\overline{x}'d\zeta\\
&\le \|x''\|^2_{L^2(0,\pi)} -\|x^{(3)}\|^2_{L^2(0,\pi)}\\
&\le 0,
\end{align*}
by the Poincaré inequality. It is easy to see that  $\ran(E - A) = \Z$ and  $\ker A \cap \ker E=\{0\}$.
Next we calculate $A^*:\D(A^*)\subset \Z\rightarrow \X$.

Note that $S:\X\rightarrow \Z$ given by $Sf:=f''$ is an  isometric isomorphism with 
$$ (S^{-1}f)(x)=\int_0^x (x-t)f(t)dt-\frac{x}{\pi}\int_0^\pi (\pi-t)f(t)dt.$$

Then, for $x\in \D(A)$ and $z\in \X$
\begin{align*}
\langle Ax, z\rangle_\Z &= \int_0^\pi x''\overline{z} +2x^{(4)}\overline{z}d\zeta\\
&= \int_0^\pi x''\overline{z} +2x{''}\overline{z}''d\zeta\\
&= \int_0^\pi x''(\overline{S^{-1}z+ 2z})''d\zeta\\
&=\langle x, A^*z\rangle_\X
\end{align*}
with $A^*z=S^{-1}z+ 2z$ for $z\in \X$. For $x\in \D(A^*)=\X$ and $y=S^{-1}x$ we calculate
\begin{align*}
\re\langle A^*x, E^*x\rangle_\X&= \re\langle EA^*x,x\rangle_\Z\\
&=\re\int_0^\pi (S^{-1}x+ x+2x+ 2x'')\overline{x}d\zeta\\
&=\re\int_0^\pi (y+ y''+2y''+ 2y^{(4)})\overline{y''}d\zeta\\
&= -\|y'\|_\Z^2+ \|y''\|_\Z^2 \\
&\quad  -2 \re\int_0^\pi y'\overline{y^{(3)}} d\zeta -2\|y^{(3)}\|_\Z^2\\
&=  \|y''\|_\Z^2 -\|y^{(3)}\|_\Z^2  \\
&\le 0.
\end{align*}
 The calculations are simpler  than using the eigenfunction expansion. Furthermore,  with this approach it can be concluded that not only is the system well-posed, but  also
 that the dynamics are a contraction with respect to the $L^2(0,\pi) $-norm. 

\section*{Conclusion}

Conditions for the solvability of partial differential-algebraic equations on Hilbert spaces have been presented and a generalization of the Lumer-Phillips  generation theorem to dissipative differential-algebraic equations.  The results were illustrated with some applications. 

It is assumed that $E$ is a linear bounded operator from $\X$ to $\Z$. In order to apply the results to some other applications such as dissipative water waves one has to deal with unbounded operators $E$. This will be the subject of future  work.

\newcommand{\etalchar}[1]{$^{#1}$}


\end{document}